\documentclass[a4paper,twoside,11pt]{amsart}
\usepackage{graphicx}
\usepackage{showlabels}
\usepackage{multicol}
\usepackage{geometry}
\usepackage{t1enc}	
\usepackage[square]{natbib} 
\usepackage[utf8]{inputenc}

\usepackage{amsmath}
\usepackage{amsfonts}
\usepackage{amssymb} 
\usepackage{color,euscript,graphicx,tabularx,epsfig,enumerate}
\usepackage{graphicx}
\usepackage{version}
\def\h1{\hspace{1cm}}
\def\h2{\hspace{2cm}}\def\h3{\hspace{3cm}}
\def\h4{\hspace{4cm}}\def\h5{\hspace{5cm}}
 \definecolor{grey}{rgb}{0.75,0.75,0.75}
\definecolor{orange}{rgb}{1.0,0.5,0.5}
\definecolor{brown}{rgb}{0.5,0.25,0.0}
\definecolor{pink}{rgb}{1.0,0.5,0.5}

\def\dis{\displaystyle}
 
\def\paragraph#1{\textit{#1}.}

\newtheorem{pp}{Proposition}

\newtheorem{lm}{Lemma}

\newtheorem{defi}{Definition}
\newtheorem{req}{Remark}

\newtheorem{pb}{Problem}

\newcommand{\faitdisparaitre}[1]{}

 \def\cqfd{\unskip\kern 6pt\penalty 500
\raise -2pt\hbox{\vrule\vbox to5pt{\hrule width 4pt
\vfill\hrule}\vrule}}

\def\A{\mathbf{A}}
\def\A2{{\mathbf{A}}^2}

\def\AB0{{\A}^2_{B( {0} ,\,{1})}}
\def\AB#1#2{{\A}^2_{B( {#1} ,\,{#2})}}

\def\R{\mathbb{R}}

\def\dim{\mathrm{dim}}

\def\monomials#1{\mathcal{M}}

\def\dis{\displaystyle}

\title[Revisiting Clifford Algebras]{Revisiting Clifford Algebras}
\author{Marc Atteia}
\address{} 
\email{marcatteia@orange.fr} 
\date{Version of \today}
\begin{document}
\maketitle 

\begin{abstract} 
This paper is divided in three parts.
In the first part, I study the Clifford algebra associated to the
hessian of a functional $f$  defined on an open subset of
 $\R^n$ \ and the Clifford algebra associated to the hessian of the Legendre
transform of $f$.
 I give also the definition of a tensorial topology on a
Clifford algebra.
 In the second part, I study the Clifford algebra $C\left( E,q\right) $
of an infinite dimensional Banach space $E$ and its main properties.
 Finally, in the third part, I give the explicit formula of the
hilbertian kernel of a Clifford algebra with examples.
\end{abstract}
\section{Basic properties of a Clifford algebra}
\subsection{Quadratic space.}
\begin{defi}
Let $\mathbb{K=R}$ or $\mathbb{C}$ \ and $E$ a
vector space on $\mathbb{K}$.
A quadratic form $q$ on \ $E$ \ is a mapping from \ $E$ \ into $\mathbb{K}$
 such that :
\begin{itemize}
\item[(i)] $\forall x\in E$ , $\forall \lambda \in \mathbb{K}$ , $q\left(
\lambda x\right) =\lambda ^{2}q\left( x\right) $
\item[(ii)] If \ $\forall x,y\in E$ , $b\left( x,y\right) =q\left(
x+y\right) -q\left( x\right) -q\left( y\right) $
then $b$ is bilinear and symmetric .
\end{itemize}
So, we say that $\left( E,q\right) $  is a quadratic space (on 
 $\mathbb{K}$ ).
 \end{defi}
\subsection{Clifford algebra }
\begin{pb}\label{P1}
Let $\left( E,q\right) $ \ be a quadratic space (on \ $\mathbb{K}$ ).
We want to find if there exists an algebra \textbf{\ A }on \ $\mathbb{K}$ ,
and
a linear mapping $\varphi $ from \ $E$ \ into $\mathbf{A}$ such that :
$$\forall x\in E , \, \left[ \varphi \left( x\right) ^{2}\right] =q\left(
x\right) .\mathbf{1}_{\mathbf{A}}$$
where $\mathbf{1}_{\mathbf{A}}$ is the identity of $\mathbf{A}$.
\end{pb}
So, if the problem P1 has a solution, then :
$$\forall x,y\in E, \, \varphi \left( x\right) .\varphi \left(
y\right) +\varphi \left( y\right) .\varphi \left( x\right) =2b\left(
x,y\right) \mathbf{1}_{\mathbf{A}}.$$

In that case, we set :
$$W=\varphi \left( E\right)\,\textrm{and} \, \forall x,y\in E,\, B\left(
\varphi \left( x\right) .\varphi \left( y\right) \right) =b\left( x,y\right)
. \mathbf{1}_{\mathbf{A}}.$$
Then, $W$ is a linear space on $\mathbb{K}$ and
 $B$ is a bilinear and symmetric mapping from $W$ into $\mathbb{
K}$ .
 \\
The problem~\ref{P1} is equivalent to the following :
\begin{pb}\label{P2}
Find an algebra $W$ \ on $\mathbb{K}$ such that :
$$\forall v,w\in W,\, v\cdot w+w\cdot v=2B\left( v,w\right) $$
where $\cdot $ is the product in $W$ .
\end{pb}
Now, let us suppose that $E$ is finite dimensionnal.
Let $n\in \mathbb{N}^{\ast }$ be the dimension of $E$ .
We can prove that :
\\
The problem~\ref{P1} has a solution and
there exists an algebra on $\mathbb{K}$ , denoted by \ $C\left( E,q\right) $
such that :
\\if $\left( e_{1},...,e_{2}\right) $ is an (orthogonal) basis of $E$ \
such that :
\\
$$\forall j,l\in \left( 1,...,n\right), \, b\left( e_{j},e_{l}\right)
=0\quad \textrm{if} \quad j\neq l\quad \textrm{and}\quad  b\left( e_{j},e_{j}\right) =\pm 1$$
then the set
$$\{\varphi \left( e_{j_{1}}\right) \cdot ...\cdot\varphi \left(
e_{j_{k}}\right)\, ;\, 1\leq j_{1}<j_{2}<...<j_{k}\leq n\}$$
is a basis of $C\left( E,q\right) $ and $\ \dim \left( C\left( E,q\right)
\right) =2^{n}$ .

\subsection{Clifford algebra and tensor analysis}
\quad\\
When $p\in \mathbb{N}$ , $p\geq 2$ , we set : $E^{\otimes p}=\underset{p}{%
\underbrace{E\otimes .......\otimes E}}$. Then
$ E^{\otimes p}$ \ is generated by the following tensors :
$$ x_{1}\otimes .......\otimes x_{p}, \  \ x_{j}\in E, \  \ 
j\in \left( 1,...,p\right). $$
 Let $E^{\wedge p}$ be the vector space generated by the following
antisymmetric tensors :
$$\frac{1}{p!}\sum_{\sigma \in \emph{G}_{p}}\varepsilon \left(
\sigma \right) \left( x_{\sigma \left( 1\right) }\otimes .......\otimes
x_{\sigma \left( p\right) }\right) $$
 where
$ \emph{G}_{p}$ is the symmetric group of order $p$ , $x_{\sigma
\left( j\right) }\in E$ , \ \ $j\in \left( 1,...,p\right) $
 and \ $\varepsilon \left( \sigma \right) =+1$ $\left( resp.-1\right) $
if $\sigma $ is even $\left( resp.\text{ odd}\right) $.
\\
So, when \ $p\in \mathbb{N}$ , $p\geq 2$ ,
 we set : $$\mathcal{E}_{p}^{\otimes }=E\oplus E^{\otimes
2}\oplus .....\oplus E^{\otimes p} \ \ \textrm{and}
\ \ \mathcal{E}_{p}^{\wedge }=E\oplus E^{\wedge
2}\oplus .....\oplus E^{\wedge p} .$$
Then, we can prove the following
\begin{pp} With the same hypotheses as in the previous paragraph,
 $C\left( E,q\right) $ is (algebra\"{\i}cally) isomorphic to \ $%
\mathcal{E}_{p}^{\wedge }$ .
\end{pp}
\subsection{Topology on $C\left( E,q\right)$}
\quad\\
Let us suppose that \ $E$ is a $n-$dimensional $\left( n\in \mathbb{N}^{\ast
}\right) $ vector space on $\mathbb{K}$ ,
and $\nu $ \ a norm on \ $E$ .
Let $\gamma $ be a tensorial (semi-) norm such that \ : \ $\pi \leq \gamma
\leq \varepsilon $ \ ( cf. $B$ $2$)
where \ $\pi $ (resp. $\varepsilon $ ) is the canonical projective (resp.
injective) tensorial (semi-)norm.
\\
 So, we set  :
$\forall p\in \mathbb{N}$ , $p\geq 2$ , \ $\nu _{\gamma }^{\otimes p}=\nu
\gamma \nu .....\gamma \nu $ \ where \ $\gamma $ \ appears \ $\left(
p-1\right) $ \ times.
\\
We denote by $\left( E^{\otimes p},\nu _{\gamma }^{\otimes p}\right) $ 
the vector space \ $E^{\otimes p}$ \ imbedded with the (semi-)norm \ $\nu_{\gamma }^{\otimes p}$ 
and \ $\mathcal{E}_{\gamma ,p}^{\otimes }$ \ the direct sum of the spaces \
\ $\left( E,\nu \right) $ , \ $\left( E^{\otimes 2},\nu _{\gamma }^{\otimes
2}\right) $ , ...., $\left( E^{\otimes p},\nu _{\gamma }^{\otimes p}\right) $.
\\
Then
 $\mathcal{E}_{\gamma ,p}^{\otimes }$ \ is a topological vector space.
\\
Let \ $\mathcal{E}_{\gamma ,p}^{\wedge }$ \ be the space \ $\mathcal{E}%
_{p}^{\wedge }$ \ imbedded with the topology induced by that of \ \ $%
\mathcal{E}_{\gamma ,p}^{\otimes }$ .
\begin{pp}
 $C\left( E,q\right) $ \ is a topological vector space when it is
imbedded by the reciprocal image
topology of that we have defined above on \ $\mathcal{E}_{\gamma ,p}^{\wedge}$ .
\end{pp}
\subsection{About second derivatives of a functional}
\begin{defi}
 Let $\ \Omega $ \ be an open set in \ $\mathbb{R}^{n}$ , $n\in 
\mathbb{N}^{\ast }$ \ and $\ f\in C^{2}\left( \Omega \right) $ .
\\
Then, if \ $a\in \Omega $ , $f^{\prime }\left( a\right) \in \mathcal{L}%
\left( \ \mathbb{R}^{n};\ \mathbb{R}\right) $ and  $f^{\prime \prime }\left( a\right) \in \mathcal{L%
}\left( \ \mathbb{R}^{n};\ \mathcal{L}\left( \ \mathbb{R}^{n};\ \mathbb{R}%
\right) \right) =\mathcal{B}\left( \mathbb{R}^{n}\times \mathbb{R}^{n};%
\mathbb{R}\right) $  which is the space of bilinear \ (and
continuous) mappings from \ $\mathbb{R}^{n}\times \mathbb{R}^{n}$ \
into \ $\mathbb{R}$ .
\\
Let $\left\langle \cdot \mid \cdot \right\rangle $  be the canonical
euclidean scalar product on \ $\mathbb{R}^{n}$ . Then 
$\forall h=\left( h_{1},....,h_{n}\right)$
and
$
\forall k=\left(
k_{1},....,k_{n}\right)\in \mathbb{R}^{n}$ we let :
 $$
f^{\prime }\left( a\right) \cdot h=\sum_{j=1}^{n}\frac{\partial f}{\partial
x_{j}}\left( a\right) \cdot h_{j}=\left\langle f^{\prime }\left( a\right)
\mid h\right\rangle $$
  and
\begin{align*}
\left( f^{\prime \prime }\left( a\right) \cdot h\right) \cdot
k&=\sum_{j,l\text{ }=1}^{n}\frac{\partial ^{2}f}{\partial x_{j\text{ }%
}\partial x_{l}}\left( a\right) .h_{j}.k_{l}
\\&=
\left\langle f^{\prime \prime
}\left( a\right) \cdot h\mid k\right\rangle 
\\&=
\left\langle f^{\prime \prime
}\left( a\right) \mid h\otimes k\right\rangle _{\otimes }.
\end{align*}
In the following, we shall identify :
 \begin{itemize}
 \item[(i)]
 $f^{\prime }\left( a\right) $ \ with the vector \ $^{t}\left( \frac{%
\partial f}{\partial x_{1}}\left( a\right) ....\frac{\partial f}{\partial
x_{n}}\left( a\right) \right) $
\\
\item[(ii)]
 $f^{\prime \prime }\left( a\right) $ \ with the matrix 
$\dis
\left[ \frac{\partial ^{2}f}{\partial x_{j\text{ }}\partial x_{l}}\left(
a\right) \right] _{j,l=1.....n}\text{ \ which is symmetric}.
$
\end{itemize}
\end{defi}
\subsubsection{Geometrical point of view}
 \quad\\
 Let \ $graph\left( f\right) =\left\{ \left(
x,z\right) \in \Omega \times \mathbb{R}\text{ ; }z=f\left( x\right) \right\} 
$
then $\left( f^{\prime }\left( a\right) ,-1\right) $ \ is a
vector in \ \ $\mathbb{R}^{n+1}$ \ which is normal
at the point $\left( a,f\left( a\right) \right) $ \ to \ $graph\left(
f\right) $ . Moreover
$$\forall y\in \Omega  \quad P_{y}=\left\{ \left( x,z\right) \in 
\mathbb{R}^{n}\times \mathbb{R}\text{ ; }z-xf^{\prime }y)-f\left( y\right)
+\left\langle y\mid f^{\prime }\left( y\right) \right\rangle =0\right\} $$
is an hyperplane which is tangent to $graph\left( f\right) $ at the point \ $%
\left( y,f\left( y\right) \right) $ .
\\
In the following, we shall denote by \ \ $q_{f}\left( \cdot ,a\right) $
the quadratic form :
$\ h\rightarrow \left\langle f^{\prime \prime }\left( a\right) \mid h\otimes
h\right\rangle _{\otimes }$.
\subsubsection{Legendre transform of a functional}
\begin{defi}
We call Legendre transform of \ \ \ $%
graph\left( f\right) $ \ the set :
$$\left( \ graph\left( f\right) \right) ^{\ast }=
\left\{ \left( x^{\ast
},z^{\ast }\right) \in \mathbb{R}^{n}\times \mathbb{R};\ x^{\ast }=f^{\prime
}y)\ \text{and}\ \ z^{\ast }=f\left( y\right) -\left\langle y\mid f^{\prime
}\left( y\right) \right\rangle \text{ },\text{ }y\in \Omega \right\} $$
\end{defi}
Generally, \ $\left( \ graph\left( f\right) \right) ^{\ast }$ \ is
not the graph of a mapping from \ $\mathbb{R}^{n}\ $\ into $\ \mathbb{R}$.
\\
In the following we set :
 $$\forall x^{\ast }\in \mathbb{R}^{n}, \quad
 f^{\ast }\left( x^{\ast
}\right) =\left\{ z^{\ast }\in \mathbb{R}\text{ ; }\ x^{\ast }=f^{\prime
}y)\ \text{and\ }\ z^{\ast }=f\left( y\right) -\left\langle y\mid f^{\prime
}\left( y\right) \right\rangle ,y\in \Omega \right\}.$$

\subsubsection{Examples}
\begin{itemize}
\item[(i)]  $\forall x\in \mathbb{R}$\ , $\ f\left( x\right)
=p^{-1}\left\vert x\right\vert ^{p}$ \ , when $p\in \mathbb{N}^{\ast }$ .
Then :
\\ $\ f^{\ast }\left( x^{\ast }\right) =-\left( p^{\ast }\right)
^{-1}\left\vert x^{\ast }\right\vert ^{p^{\ast }}$ \ with \ $p^{\ast }\in 
\mathbb{N}^{\ast }$ \ and \ $p^{-1}+\left( p^{\ast }\right) ^{-1}=1$ , \ $x^{\ast }\in \mathbb{R}$ .
\item[(ii)] $\forall x\in \mathbb{R}$\ , $\ f\left( x\right)
=\left( x^{2}-1\right) ^{2}$ . Then :
\\ $\left( \ graph\left( f\right) \right) ^{\ast }=$\ $\left\{ \left(
x^{\ast },z^{\ast }\right) \in \mathbb{R}\times \mathbb{R};\ x^{\ast
}=4y\left( y^{2}-1\right) \ \text{and}\ \ z^{\ast }=\left( y^{2}-1\right)
\left( 3y^{2}+1\right) \text{ },\text{ }y\in \mathbb{R}\right\} $
\item[( iii)] Minkowski quadratic form.
\\ Let \ $\forall x=\left( x_{1}....x_{n}\right) $, $f\left( x\right) =\left( x_{1}\right) ^{2}+....+\left(
x_{p}\right) ^{2}-\left( \left( x_{p+1}\right) ^{2}+....+\left( x_{n}\right)
^{2}\right) $ \ when \ $1\leq p\leq n$ .
Then :
$\forall a\in \mathbb{R}^{n}$ ,
$f^{\prime \prime }\left( a\right) $ \ is a diagonal matrix
with elements are equal to $\left( +1\right) $ \ $\left( resp.\text{ }\
-1\right) $
and $$\forall a^{\ast }\in \mathbb{R}^{n},\quad \left( f^{\ast
}\right) ^{\prime \prime }\left( a^{\ast }\right) =\frac{1}{2}f^{\prime
\prime }\left( a\right).$$
\end{itemize}
\subsubsection{The relation between \ $f^{\prime \prime }$ \
and \ $\left( f^{\ast }\right) ^{\prime \prime }$ }
\quad\\
 Let us suppose that $f^{\prime }$ \ is a diffeomorphism from a
(non void) open subset $\Omega _{0}$
contained in \ $\Omega $ \ \ onto an other (non void) open subset $\Omega_{0}^{\ast }\subset $ $\mathbb{R}^{n}$ .
 With previous notations, we have :
$$\left( x^{\ast }=f^{\prime }(y)\text{ , }y\in \Omega
_{0}\right) \Longrightarrow \left( y=\left( f^{\prime }\right) ^{-1}\left(
x^{\ast }\right) \underset{def}{=}g\left( x^{\ast }\right) \text{ ; }x^{\ast
}\in \Omega _{0}^{\ast }\right). $$
 Then :
$$z^{\ast }=f\left( g\left( x^{\ast }\right) \right)
-\left\langle g\left( x^{\ast }\right) \mid f^{\prime }\left( g\left(
x^{\ast }\right) \right) \right\rangle =f^{\prime \prime }\left( x^{\ast
}\right)\,;\, x^{\ast }\in \Omega _{0}^{\ast }.$$
 It can be proved easily that  :
 \\
when \ $n\in \mathbb{N}$ , $n\geq 2$ and when $\left( f\right) ^{\prime
\prime }\left( y\right) $ \ is an invertible operator, then :
$$\left( f^{\ast }\right) ^{\prime \prime }\left( x^{\ast }\right)
=\left( \left( f\right) ^{\prime \prime }\left( y\right) \right) ^{-1},\ \ 
y\in \Omega _{0},\ \ x^{\ast }\in \Omega _{0}^{\ast }.$$

\subsection{ Clifford algebra associated to the second
derivative of a functional}
\quad \\
Let \ $f\in C^{2}\left( \Omega \right) $ .
For all $a\in \Omega $ , \ we denote by \ $B_{f}\left( \cdot ,\cdot \text{ ; 
}a\right) $ \ the bilinear form which is associated to
$q_{f}\left( \cdot ,a\right) $ and by \ $M_{f}\left( a\right) $ \ the $%
n\times n$ \ matrix 
$ \dis \left[ \frac{\partial ^{2}f}{\partial x_{j\text{ }}\partial x_{l}}\left(
a\right) \right] _{j,l=1.....n}
$.
 Let us suppose that \ $M_{f}\left( a\right) $ \ is invertible.
As \ $M_{f}\left( a\right) $ \ is a symmetric matrix, then \ $M_{f}\left(
a\right) $ \ has $\ n$ \ eigenvectors $e_{1},...,$ $e_{n}$
which are linearly independant such that :
\\ $B_{f}\left( e_{j},e_{l}\text{ ; }a\right) =0$ \ if \ $j\neq l$
\ , \ $1\leq j,l\leq n$ and \ $B_{f}\left( e_{j},e_{j}\text{ ; }a\right) \neq 0$ \ \ \
if \ $1\leq j\leq n$ . 
So,
\begin{pp}
To the functional \ $f$ \ defined above, we can associate the Clifford
algebras :
\\
\qquad $C\left( \mathbb{R}^{n},q_{f}\left( \cdot ,a\right) \right) $ \ , \ $%
a\in \Omega _{0}$ \ \ and \ \ $C\left( \mathbb{R}^{n},q_{f^{\ast }}\left(
\cdot ,a^{\ast }\right) \right) $ \ , \ $a^{\ast }\in \Omega _{0}^{\ast }$ .
\end{pp}
\pagebreak
\section{Clifford algebra on Banach
spaces}
\subsection{ Definitions}
\subsubsection{ Banach spaces}
\quad\\
 Let \ $B$ \ be a linear space and \ $\nu $ \ a norm on \ $B$ .
\\
We say that $\left( B,\nu \right) $ \ is a normed space.
\\We denote by \ $T\left( B,\nu \right) $ \ the vectorial topology on $B$
associated to \ $\nu $ .
\\
We say that \ $\left( B,\nu \right) $ \ is a Banach space when \ $\left(
B,\nu \right) $ with $T\left( B,\nu \right) $	
is complete.
\\
 We shall denote by \ \ $B^{\ast }$ the topological dual of \ $\left(
B,\nu \right) $
and by \ $\nu ^{\ast }$ the dual norm of \ $\nu $ \ on \ $B^{\ast }$ .
\\
 We denote, also, by \ $\left\langle \cdot ,\cdot \right\rangle $ \
the duality bracket between \ $B$ \ and $B^{\ast }$ .

\subsubsection{Schauder basis}
\quad\\
 Let \ $\left( B,\nu \right) $ \ be an infinite dimensional Banach
space.
We say that a sequence $\left\{ e_{j}\text{ ; }j\in \mathbb{N}^{\ast
}\right\} \subset B$ \ is a Schauder basis in \ $\left( B,\nu \right) $ if 
for any $\ x\in B$\ , there exists a unique sequence \ $\left\{
\alpha _{j}\left( x\right) \text{ ; }j\in \mathbb{N}^{\ast }\right\} \subset 
\mathbb{K}$
such that : \ $\underset{n\rightarrow \infty }{\lim }$ $\nu \left(
x-\sum_{j=1}^{n}\alpha _{j}\left( x\right) e_{j}\right) =0$ .
\\
 Below, we use the following notations :
$\forall x\in B$ \ , \ $\forall n\in \mathbb{N}^{\ast }$ , \ 
\\ $\beta _{n}\left( x\right) =\sum_{j=1}^{n}\alpha _{j}\left(
x\right) e_{j}$ \ , \ $\dis \rho _{n}\left( x\right) =\sum_{j=n+1}^{\infty
}\alpha _{j}\left( x\right) e_{j}$ \ and \ $B_{n}=\left\{ \beta _{n}\left(
x\right) \text{ \ ; }x\in B\right\} $ .
\\
Then $B_{n}$ \ is a \ linear subspace in $B$ .
 Moreover, we suppose that the topology on $B_{n}$ \ is induced by the
topology $T\left( B,\nu \right) $ .
\\
It is clear that :

\qquad $\forall n\in \mathbb{N}^{\ast }$ , \ $B_{n}\subset B_{n+1}\subset B$
\ and \ 

$\qquad \underset{n}{\cup}\left( B_{n}\right) =B$ \ (the sequence \ $\left(
B_{n}\right) $ pointwise converges to $B$).

\qquad The injection of \ $B_{n}$ \ into \ $B_{n+1}$ (or in \ $B$ ) \ is
continuous.
\subsection{ Tensorial topological products of Banach spaces}
\subsubsection{Definitions}
\quad\\
 Let $\left( E,\lambda \right) $ \ and \ $\left( F,\mu \right) $ \ be
two (semi-) normed linear spaces.
 The \textit{projective }tensorial product of \ $\left(
E,\lambda \right) $ \ and \ $\left( F,\mu \right) $ \ is the (semi-)normed
linear space \ $\left( \ E\otimes F,\left( \lambda \pi \mu \right)
\right) $ \ where \ $\lambda \pi \mu $ \ is the tensorial (semi-) norm on \ $E\otimes F$
such that :
$$\forall z\in E\otimes F, \quad\left( \lambda \pi \mu
\right) \left( z\right) =Inf\left\{ \sum_{k}\lambda \left( x_{k}\right) \mu
\left( y_{k}\right) \text{ ; \ }z=\sum_{k}x_{k}\otimes y_{k}\right\}. $$
 Let $E^{\ast }$ \ \ $\left( resp.\text{ \ }F^{\ast }\right) $ \ be
the topological dual of \ $\left( E,\lambda \right) $ \ $\left( resp.\text{ }
\left( F,\mu \right) \right) $ \ and
$\lambda ^{\ast }$ \ $\left( resp.\text{ \ }\mu ^{\ast }\right) $ \
be the dual norm of \ $\lambda $ \ $\left( resp.\text{ }\mu \right) $ . Then :
\\ the \textit{inductive} tensorial product of \ $\left(
E,\lambda \right) $ \ and \ $\left( F,\mu \right) $ \ is the (semi-) normed linear space \ $\left( \ E\otimes F,\left( \lambda \varepsilon \mu
\right) \right) $ \ where \ $\lambda \varepsilon \mu $ \ is the tensorial
(semi-) norm on \ $E\otimes F$ such that :
\\ $$\forall z\in E\otimes F, \quad
\left( \lambda \varepsilon \mu
\right) \left( z\right) =Sup\left\{ \left\vert \sum_{k}\left\langle
x_{k},x^{\ast }\right\rangle \left\langle y_{k},y^{\ast }\right\rangle
\right\vert \text{ ; }z=\sum_{k}x_{k}\otimes y_{k}\right\} $$
 where\ \ $x^{\ast }\in E^{\ast }$ \ , \ \ $y^{\ast }\in F^{\ast }$ , $%
\lambda ^{\ast }\left( x^{\ast }\right) \leq 1$ \ , \ $\mu ^{\ast }\left(
y^{\ast }\right) \leq 1$.
\\ It can be proved that : \ $\varepsilon \leq \pi $ , and we say
that a tensorial norm \ $\gamma $ \ on $E\otimes F$
is a \textit{reasonable} tensor norm(r.t.n) \ if : \ $\varepsilon
\leq \gamma \leq \pi $ .
\subsection{Schauder bases on topological tensorial products}
\quad\\
 Let $\left( e_{j}\right) $ \ $\left( resp.\text{\ }\left(
f_{k}\right) \right) $ a Schauder basis of \ $\left( E,\lambda \right) $ \ $%
\left( resp.\text{ }\left( F,\mu \right) \right) $ .
Then, if \ $x\in E$ \ $\left( resp.\text{\ }y\in F\right) $
$\dis
x=\sum_{j=1}^{\infty }\zeta _{j}\left( x\right) e_{j}$ \ $\dis
\left(
resp.\text{ }\ y=\sum_{k=1}^{\infty }\eta _{k}\left( y\right) f_{k}\text{ }
\right) $
we set : 
\\
$\forall n\in \mathbb{N}^{\ast }$ ,
 $\forall x\in E$ \ , 
 $\dis\varphi _{n}\left( x\right)
=\sum_{j=1}^{n}\zeta _{j}\left( x\right) e_{j}$ \ and \ \ $\dis
\chi _{n}\left(
x\right) =\sum_{j=n+1}^{\infty }\zeta _{j}\left( x\right) e_{j}$.
\\
 $\dis\left( resp.\text{ }\forall y\in F,\text{ \ }\psi _{n}\left(
y\right) =\sum_{k=1}^{n}\eta _{k}\left( y\right) f_{k}\ \ \ \text{and }\ \
\omega _{n}\left( y\right) =\sum_{kj=n+1}^{\infty }\eta _{k}\left( y\right)
f_{k}\right) $.
\\
So, we can easily prove that : 
$$\underset{n\rightarrow \infty }{\lim }\lambda \left( \chi
_{n}\left( x\right) \right) =\underset{n\rightarrow \infty }{\lim } 
\mu \left( \omega _{n}\left( y\right) \right) =0.$$

\subsubsection{First case: \ $E\otimes F$ }
 \quad\\
$\forall n\in \mathbb{N}^{\ast }$ , $\forall x\in E$ \ , \ $%
\forall y\in F$,
 $$x\otimes y=\left( \varphi _{n}\left( x\right) +\chi _{n}\left(
x\right) \right) \otimes \left( \text{\ }\psi _{n}\left( y\right) +\omega
_{n}\left( y\right) \right) =A_{n}\left( x,y\right) +R_{n}\left( x,y\right) $$
where \ : \ $A_{n}\left( x,y\right) =\varphi _{n}\left(
x\right) \otimes \psi _{n}\left( y\right) $ \ \ and \ 

\qquad\quad  $R_{n}\left( x,y\right) =\left( \varphi
_{n}\left( x\right) \otimes \omega _{n}\left( y\right) \right) +\left( \chi
_{n}\left( x\right) \otimes \text{\ }\psi _{n}\left( y\right) \right)
+\left( \left( \chi _{n}\left( x\right) \otimes \text{\ }\omega _{n}\left(
y\right) \right) \right) $
\\
Then, if \ $\gamma $ \ is a \textit{reasonable} tensor norm on \ $%
E\otimes F$ , we have :
$$\gamma \left( R_{n}\left( x,y\right) \right) \leq \lambda
\left( \varphi _{n}\left( x\right) \right) \mu \left( \omega _{n}\left(
y\right) \right) +\lambda \left( \chi _{n}\left( x\right) \right) \mu \left(
\psi _{n}\left( y\right) \right) +\lambda \left( \chi _{n}\left( x\right)
\right) \mu \left( \omega _{n}\left( y\right) \right) $$
and
$$\underset{n\rightarrow \infty }{\lim }
\left(
x\otimes y\right) -A_{n}\left( x,y\right) =0 .$$
 We verify easily that : \ $\dis A_{n}\left( x,y\right)
=\sum_{j,k=1}^{n}\zeta _{j}\left( x\right) \eta _{k}\left( y\right) \left(
e_{j}\otimes f_{k}\right) $ .
\\
Now, we consider the infinite sequence \ $J=\left\{ \left(
m,n\right) \text{ ; }m,n\in \mathbb{N}^{\ast }\right\} $.
We set :
\begin{align*}
J_{1}&=\left( 1,1\right) 
\\
J_{2}&=\left( 1,2\right) \left( 2,1\right) \left( 2,2\right) 
 \\
 J_{3}&=\left( 1,3\right) \left( 2,3\right) \left( 3,3\right)
\left( 3,2\right) \left( 3,1\right) 
\\
&\vdots
\\
J_{l}&=\left( 1,l\right) \left( 2,l\right) \left( 3,l\right)
.....\left( l,l\right) \left( l,l-1\right) ....\left( l,2\right) \left(
l,1\right),\quad l\in \mathbb{N}^{\ast }.
\end{align*}
Then : $J=\underset{l\in \mathbb{N^\ast}}{\cup} J_{l}$ .
\begin{pp}
The sequence $\left\{ \left( e_{j}\otimes f_{k}\right) \text{
, \ }\left( j,k\right) \in J_{l}\text{ \ , }l\in \mathbb{N}^{\ast }\text{\ }%
\right\} $ \ is a Schauder basis of the normed
$\left( E\otimes F\text{ \ , \ }\lambda \gamma \mu
\right) $ .
 \end{pp}
 \subsubsection{Second case : \ $E\wedge F$ }
\quad\\
 Let \ $E\wedge F$ \ be the linear space generated by the
vectors such that 
$$x\wedge y=\frac{1}{2!}\left( x\otimes y-y\otimes
x\right) $$
Using previous notations, we have :
 $2\left( x\wedge y\right) =A_{n}^{\wedge }$ $\left( x,y\right)
+R_{n}^{\wedge }$ $\left( x,y\right) $ \ with :
 \begin{align*}
 A_{n}^{\wedge }\left( x,y\right) &=\left(
\sum_{j=1}^{n}\zeta _{j}\left( x\right) e_{j}\right) \otimes \left(
\sum_{k=1}^{n}\eta _{k}\left( y\right) f_{k}\ \right) -\left(
\sum_{k=1}^{n}\eta _{k}\left( y\right) f_{k}\ \right) \otimes \left(
\sum_{j=1}^{n}\zeta _{j}\left( x\right) e_{j}\right) 
\\
&=\sum_{j,k=1}^{n}\zeta _{j}\left(
x\right) \eta _{k}\left( y\right) \left( e_{j}\wedge f_{k}\right) 
\end{align*}
 and  $$\underset{n\rightarrow \infty }{\lim \text{ }}
R_{n}^{\wedge }\left( x,y\right) =0.$$
 As in the first case, we can prove the
\begin{pp}
The sequence $\left\{ \left( e_{j}\wedge f_{k}\right) \text{ ,
\ }\left( j,k\right) \in J_{l}\text{ \ , }l\in \mathbb{N}^{\ast }\text{\ }%
\right\} $ \ is a Schauder basis of the normed
space \ \ $\left( E\wedge F\text{ \ , \ }\lambda \gamma \mu
\right) $ .
\end{pp}
\subsubsection{Third \ case : \ $E\vee F$}
\quad\\
 Let \ $E\vee F$ \ be the linear space generated by the vectors
such that $$x\vee y=\frac{1}{2!}\left( x\otimes y+y\otimes
x\right). $$
As in the second case, we can prove the
\begin{pp}
 The sequence $\left\{ \left( e_{j}\vee f_{k}\right) \text{ , \ 
}\left( j,k\right) \in J_{l}\text{ \ , }l\in \mathbb{N}^{\ast }\text{\ }
\right\} $ \ is a Schauder basis of the normed
 space \ \ $\left( E\vee F\text{ \ , \ }\lambda \gamma \mu
\right) $ .
\end{pp}
 \subsubsection{Extension. Exercice.}
\qquad\\ Let $E_{j}$ be a Banach space with the Schauder basis \ $%
\left( e_{j}^{k}\right) _{k\in \mathbb{N}^{\ast }}$ \ $j=1....n$ .
Deduce from the previous paragraphs what is the Schauder bases of
the spaces  $E_{1}\otimes \ldots\otimes E_{m}$, \ \ $E_{1}\wedge
\ldots\wedge E_{m}$ \ and \ \ \ $E_{1}\vee \ldots\vee E_{m}$ .
\subsection{Fock spaces and Schauder bases}
\subsubsection{Definitions}
\begin{itemize}
\item[(i)] Let \ $E$ \ a linear space on $\mathbb{\ K}$ .
The space 
$$\mathcal{F}^{\otimes }\left( E\right) =\mathbb{K\times }E 
\mathbb{\times }E^{\otimes 2}\times E^{\otimes 3}\times .....\times  
E^{\otimes p}\times ....$$
is called a Fock space.
 This space is an algebra.
\item[(ii)] Let \ $p\in \mathbb{N}$ , $p\geq 2$ . We set :
$$\forall x_{1},.....,x_{p}\in E, \quad 
x_{1}\vee
\ldots
\vee x_{p}=\frac{1}{p!}\sum_{\sigma \in G_{p}}\left( x_{\sigma \left(
1\right) }\otimes\ldots 
\otimes x_{\sigma \left( p\right) }\right) $$
  and
$$x_{1}\wedge
\ldots
\wedge x_{p}=\frac{1}{p!}\sum_{\sigma \in G_{p}}\varepsilon \left(
\sigma \right) \left( x_{\sigma \left( 1\right) }\otimes .....\otimes
x_{\sigma \left( p\right) }\right) $$ 
 where \ $G_{p}$ \ is the symmetric group of order \ $p$ \ and
\ $\varepsilon \left( \sigma \right) $ \ is equal to
 $\left( +1\right) $ \ $\left( resp.\text{ }\left( -1\right)
\right) $ \ if \ $G_{p}$ $\ $is even $\ \left( resp.\text{ odd}\right) $ .
\\
 We shall denote by \ $E^{\vee p}$ \ $\left( resp.\text{ }E^{\wedge p}%
\text{\ }\right) $ the linear subspace of \ $E^{\otimes p}$

\qquad which is generated by the vectors of type \ \ $x_{1}\vee \ldots\vee
x_{p}$ \ $\left( resp.\text{ }x_{1}\wedge \ldots
\wedge x_{p}\right) $ .
 So, we set :
 $$\mathcal{F}^{\vee }\left( E\right) =\mathbb{K\times }E
 \mathbb{\times }
 E^{\vee 2}\times E^{\vee 3}\times\ldots
 \times 
E^{\vee p}\times \ldots$$
 and
$$\mathcal{F}^{\wedge }\left( E\right) =\mathbb{K\times }
E
\mathbb{\times }
E^{\wedge 2}\times 
E^{\wedge 3}\times
\ldots\times E^{\wedge p}\times\ldots $$
\end{itemize}
\subsubsection{Topologies}
\qquad\\
Let $\ \left( E_{j}\text{ , }\lambda _{j}\right) $
\ $j=1,2,3$ \ three Banach spaces (on $\mathbb{K}$) \ and \ $\gamma $ \ be a
r.t.n.
 It can be proved (easily) that :
$$\left( \lambda _{1}\gamma \lambda _{2}\right) \gamma \lambda
_{3}=\lambda _{1}\gamma \left( \lambda _{2}\gamma \lambda _{3}\right) 
\underset{def}{=}\lambda _{1}\gamma \lambda _{2}\gamma \lambda _{3}.$$
Now, let $\left( E\text{ , }\lambda \right) $ \ be a Banach space (on 
$\mathbb{K}$) .
\\
 We set :
 $\forall j\in \mathbb{N}$ , \ $j\geq 2$ , $\lambda ^{\otimes
\left( j,\gamma \right) }=\lambda \gamma .....\gamma \lambda $ \ where $%
\gamma $ \ appears $\left( j-1\right) $ \ times .
 \\
 We denote by \ $\left( E\text{ , }\lambda \right) ^{\otimes j}$
(resp. $\left( E\text{ , }\lambda \right) ^{\vee j}$ , \ $\left( E\text{ , }%
\lambda \right) ^{\wedge j}$ ) the Banach space

\qquad $\left( E^{^{\otimes j}}\text{ , }\lambda ^{\otimes\left( j,\gamma \right)
}\right) $ \ ( resp. \ $\left( E^{^{\vee j}}\text{ , }\lambda ^{\otimes
\left( j,\gamma \right) }\right) $ \ , \ $\left( E^{^{\wedge j}}\text{ , }%
\lambda ^{\otimes \left( j,\gamma \right) }\right) $ .
\\ With above hypotheses, we deduce from the results in paragraph 2 the
\begin{pp}
 Let \ $\left( E,\lambda \right) $ \ be a normed space (on $\mathbb{K}$%
) \ with a Schauder basis. Then~:
 $\left( E\text{ , }\lambda \right) ^{\otimes j}$ (resp. $%
\left( E\text{ , }\lambda \right) ^{\vee j}$ , \ $\left( E\text{ , }\lambda
\right) ^{\wedge j}$ ) is a normed space with
 a Schauder basis.
\end{pp}
 Let us denote by \ $\bot $ \ one of the operators \ $\otimes $ , $%
\vee $ , $\wedge $ . Below, we set :
 $$\mathcal{F}_{\gamma }^{\bot }\left( E,\lambda \right) =\mathbb{%
K\times }\left( E,\lambda \right) \mathbb{\times }\left( E,\lambda
\right) ^{\bot 2}\times \left( E,\lambda \right) ^{\bot 3}\times 
\ldots \times 
\left( E,\lambda \right) ^{\bot j}
\times \ldots$$
 We shall say that :
 $\mathcal{F}_{\gamma }^{\bot }\left( E,\lambda \right) $ \ \
is a (tensorial) Fock space associated to \ $\left( E,\lambda \right) $
with the r.t.n. $\gamma $ .

\subsection{Limit of Fock spaces associated to a Banach space with a
Schauder basis}
\qquad\\ Let $\left( B,\nu \right) $ \ be a Banach space with a Schauder basis
\ $\left\{ e_{j}\text{ ; }j\in \mathbb{N}^{\ast }\right\} $ .
We know (cf. paragraph 1) that :

\qquad \qquad $\forall n\in \mathbb{N}^{\ast }$ , \ $B_{n}\subset
B_{n+1}\subset B$ \ and that : \ $B=\cup _{n}B_{n}$ .
\\ From the three propositions stated in the previous paragraph,  we
deduce
(easily) that : 
$$\forall n\in \mathbb{N}^{\ast }, \ \mathcal{F}_{\gamma
}^{\bot }\left( B_{n},\nu \right) \subset \mathcal{F}_{\gamma }^{\bot
}\left( B_{n+1},\nu \right) \subset \mathcal{F}_{\gamma }^{\bot }\left(
B,\nu \right) $$
and that :
$$ \mathcal{F}_{\gamma }^{\bot }\left( B,\nu
\right) =\cup _{n}\mathcal{F}_{\gamma }^{\bot }\left( B_{n},\nu \right) .$$
Now, let $q$  be a quadratic form on  $B$ \ and $q_{n}$ its
restriction to \ $B_{n}$.  We have :

$\forall x\in B_{n}$ \ , \ $q_{n}\left( x\right) =q\left(
x\right) $ . Moreover,
\\ $\forall x\in B$ \ , \ $\exists $ $n_{0}\in \mathbb{N}^{\ast }$ such
that : \ $\forall n\geq n_{0}$ \ , \ $x\in B_{n}$ \ \ \ and \ $q_{n}\left(
x\right) =q\left( x\right) $ .
\\
 Thus : $\underset{n}{lim}\,q_{n}=q$ . So,
\begin{defi}
We call Clifford algebra associated to \ $\left( B\text{ , }q\right) $
 the limit when \ $n$ \ tends to infinity
 of the Fock space \ $\left( C\left( \ B_{n}\text{ , }q_{n}\right)
\right) $ \ which can be identify to \ $\mathcal{F}_{\gamma }^{\wedge
}\left( \left( B\text{ , }\nu \right) \right) $\ with the
quadratic form $q$ \ on \ $B$ .
\end{defi}
 The results of the paragraph 1.6 can be extended (easily) to the Hessian of a regular convex functionnal.
\section{ Hilbertian Clifford algebras.}
\subsection{ Definitions.}
\subsubsection{Hilbert spaces}
\qquad \\Let \ $\mathbb{K=R}$ \ or \ $\mathbb{C}$ .
 A couple \ $\left( E,\left\langle \cdot \mid \cdot \right\rangle
\right) $ \ is a (pre-)\textit{Hilbert space} if \ $E$ \ is a linear space
on \ $\mathbb{K}$ \ and
$\left\langle \cdot \mid \cdot \right\rangle $ \ is a scalar product
on \ $E\times E$ \ such that :
\begin{itemize}
\item[(i)] $\forall y\in E$ \ \ $\left( resp.\forall
x\in E\right) $ \ , \ the mapping \ \ $\left\langle \cdot \mid
y\right\rangle $ \ \ $\left( resp.\text{ }\left\langle x\mid \cdot
\right\rangle \right) $ from \ $E$ \ into \ \ $\mathbb{K}$ \ is linear \ (resp. linear if 
$\mathbb{K=R}$ , antilinear if $\mathbb{K=C}$ ).
\item[(ii)] $\forall x\in E$ \ , \ $\left\langle x\mid
x\right\rangle \in \mathbb{R}_{+}$
\item[(iii)] If \ $x\in E$ \ \ and \ $\left\langle
x\mid x\right\rangle =0$ \ , then : $x=0$ .
\end{itemize}
 So, it can be proved that \ $\sqrt{\left\langle x\mid x\right\rangle }
$ \ is a \textit{norm} on $E$ \ denoted by \ $\left\Vert \cdot \right\Vert $.
\\
 Below, we write equivalently : \ $\left( E,\left\langle \cdot \mid
\cdot \right\rangle \right) $ \ or \ $\left( E,\left\Vert \cdot \right\Vert
\right) $ .
 If \ $\left( E,\left\Vert \cdot \right\Vert \right) $ \ is
complete, then \ $\left( E,\left\langle \cdot \mid \cdot \right\rangle
\right) $ \ is a Hilbert space.
 Thus a Hilbert space is a Banach space.
\subsubsection{Basis of a Hilbert space}
 \qquad \\ Let \ $\left( E,\left\langle \cdot \mid \cdot \right\rangle
\right) $ \ be a Hilbert space on $\mathbb{K}$ \ and \ $J=\left\{
1,....,n\right\} $ , \ $n\in \mathbb{N}^{\ast }$ \ or \ $J=\mathbb{N}^{\ast
} $\ .
 A family of elements of \ $E$ ,\ $\left\{ e_{j}\text{ ; }j\in
J\right\} $ \ is called a \textit{basis} of \ $\left( E,\left\langle \cdot
\mid \cdot \right\rangle \right) $
if for any \ $\ x\in E$ \ there exists a unique family of scalar \ $%
\left\{ \alpha _{j}\left( x\right) \text{ ; }j\in J\right\} $
such that :
$$x=\sum_{j=1}^{n}\alpha _{j}\left( x\right) e_{j}$$
when  $J=\left\{ 1,....,n\right\} $ \ and
$$\underset{k\rightarrow \infty }{\lim }\left\Vert
x-\sum_{j=1}^{k}\alpha _{j}\left( x\right) e_{j}\right\Vert $$ when 
 $J=\mathbb{N}^{\ast }$.
\\
 We say that \ $\left\{ e_{j}\text{ ; }j\in J\right\} $ \ is an 
\textit{orthonormal basis} of \ $\left( E,\left\langle \cdot \mid \cdot
\right\rangle \right) $ \ if
$\left\langle e_{j}\mid e_{l}\right\rangle =\delta _{jl}$ \ ,
\ $j,l\in J$ \ .
\begin{req} 
As any basis in an infinite dimensional Hilbert space is a Schauder
basis in that space,
 all the results stated in the part B, for Banach spaces are true for
Hilbert spaces.
\end{req}
\subsection{Reproducing Kernel and Schwartz Kernel of a Hilbert space.}
\subsubsection{Riesz theorem}
\qquad\\ Let \ $\left( E,\left\langle \cdot \mid \cdot \right\rangle \right) $
\ be a Hilbert space on $\mathbb{K}$ \ and \ $E^{\ast }=\mathcal{L}\left(
\left( E,\left\Vert \cdot \right\Vert \right) \text{ , }\mathbb{K}\right) $
its topological dual. Then :

\qquad $\forall u\in E^{\ast }$ , \ $\exists $ (only one) \ $x^{\ast }\left(
u\right) $ such that : \ $\forall x\in E$ \ , \ $u\left( x\right)
=\left\langle x\mid x^{\ast }\left( u\right) \right\rangle $ .
\subsection{Reproducing kernel or Aronszjan kernel}
\qquad \\
Let \ $\Omega $ \ be a (non void) set.
	 We say that the hilbertian space \ $\left( E,\left\langle \cdot \mid
\cdot \right\rangle \right) $ \ is a hilbertian subspace of \ $\mathbb{K}%
^{\Omega }$
and we set \ $\left( E,\left\langle \cdot \mid \cdot \right\rangle \right)
\in Hilb\left( \mathbb{K}^{\Omega }\right) $ \ if :
\begin{itemize}
\item[(i)] \ $E$ \ is a linear subspace of \ \ $\mathbb{K}%
^{\Omega }$
\item[(ii)] $\forall $ $t\in \Omega $ , $\ \exists $ $%
M_{t}\in \mathbb{R}_{+}^{\ast }$ \ such that \ $\forall x\in E$ \ , $%
\left\vert x\left( t\right) \right\vert \leq M_{t}\left\Vert x\right\Vert $
\end{itemize}
\begin{pp}
 Let $\left( E,\left\langle \cdot \mid \cdot \right\rangle \right)
\in Hilb\left( \mathbb{K}^{\Omega }\right) $. Then :
\\
 $\forall $ $t\in \Omega $ , $\exists $ (only one) \ $\mathcal{E%
}\left( \cdot ,t\right) \in E$ \ \ such that : \ $\forall x\in E$ \ , $%
x\left( t\right) =\left\langle x\mid \mathcal{E}\left( \cdot ,t\right)
\right\rangle $ .
\\
 The mapping from \ $\Omega \times \Omega $ \ into \ $\mathbb{K}$ \ :
\ $\left( s,t\right) \rightarrow \mathcal{E}\left( s,t\right) $ \ is called
the \textit{reproducing
kernel or the Aronszjan Kernel of} $\left( E,\left\langle \cdot
\mid \cdot \right\rangle \right) $.
\end{pp}
\subsubsection{Schwartz kernel}
\qquad\\
 Let \ $\mathcal{A}$ \ be a lcs (locally convex separable space) on \ $%
\mathbb{K}$ ,
$\mathcal{A}^{\ast }$ \ its topological dual and \ \ $\left\langle \cdot 
\text{ },\cdot \right\rangle $ \ the duality bracket between \ \ $\mathcal{A}
$ \ and $\mathcal{A}^{\ast }$ .
\\
 We say that the hilbertian space \ $\left( E,\left\langle \cdot \mid
\cdot \right\rangle \right) $ \ is a hilbertian subspace of \ $\mathcal{A}$
and we set \ $\left( E,\left\langle \cdot \mid \cdot \right\rangle \right)
\in Hilb\left( \mathcal{A}\right) $ \ if :
\begin{itemize}
\item[(i)]  \ $E$ \ is a linear subspace of \ $\mathcal{A}$

\item[(ii)] \ $\forall $ $a^{\ast }\in \mathcal{A}^{\ast }$ , 
$\ \exists $ $M\left( a^{\ast }\right) \in \mathbb{R}_{+}^{\ast }$ \ such
that \ $\forall x\in E$ \ , $\left\vert \left\langle jx\text{ },a^{\ast
}\right\rangle \right\vert \leq M\left( a^{\ast }\right) \left\Vert
x\right\Vert $
 where \ $j$ \ is the injection of \ $E$ \ into \ \ $\mathcal{A}
$ .
\end{itemize}
\qquad Let \ $\Lambda $ \ be the duality mapping of \ $\left( E,\left\langle
\cdot \mid \cdot \right\rangle \right) $ . Then :
$$\mathcal{A}^{\ast }\underset{j^{\ast }}{\rightarrow }%
E^{\ast }\underset{\Lambda }{\rightarrow }E\underset{j}{\rightarrow }%
\mathcal{A}$$
 where $\ \ j^{\ast }$ is the transpose of \ $j$ \ if \ $\mathbb{K=R}$
and $\ j^{\ast }$ is the conjugate transpose \ of \ $j$ \
if \ $\mathbb{K=C}$ .
\\
 So,
 $\mathcal{E=}\ j^{\ast }\Lambda \ j$ \ is called the \textit{%
Schwartz (hilbertian) Kernel \ of }$\left( E,\left\langle \cdot \mid \cdot
\right\rangle \right) $ .
\\
$\mathcal{E}$ \ is characterized by the two following
properties :
\begin{itemize}
\item[(i)] \ $\forall $ $a^{\ast },b^{\ast }\in 
\mathcal{A}^{\ast }$ \ , \ \ $\left\langle \mathit{\ }\mathcal{E}b^{\ast }
\text{ },a^{\ast }\right\rangle =
\left\{
\begin{array}{cc}
\left\langle \mathit{\ }\mathcal{E}a^{\ast }
\text{ },b^{\ast }\right\rangle &\textrm{if \ \ $\mathbb{K=R}$}
\\\\
\overline{
\left\langle \ \mathcal{E}a^{\ast }\text{ },b^{\ast }\right\rangle }
&\textrm{ if \ \ $\mathbb{K=C}$ .}
\end{array}\right .
$
\item[(ii)] $\forall $ $a^{\ast }\in \mathcal{%
A}^{\ast }$ , \ $\left\langle \mathit{\ }\mathcal{E}a^{\ast }\text{ }%
,a^{\ast }\right\rangle \geq 0$ .
\end{itemize}
\subsection{Examples of hilbertian kernels}
\subsubsection{Example 1}
\qquad\\
 Let \ $a,b\in \mathbb{R}$ \ , \ $-\infty <a<b<+\infty $ \ and \ $c\in
\left( a,b\right) $ .
We denote by $\mathcal{H}$ \ the linear space of polynomials of one variable
with degree
less or equal to \ $n\in \mathbb{N}^{\ast }$ .
\\
 We set : \ $\forall \omega _{1,}\omega _{2}\in \mathcal{H}$ , \ \ $%
\left\langle \omega _{1}\mid \omega _{2}\right\rangle =\sum_{j=0}^{n}\omega
_{1}^{\left( j\right) }\left( c\right) .\omega _{2}^{\left( j\right) }\left(
c\right) $ . Then :
$$\left( \text{\ }\mathcal{H}\text{ },\left\langle \cdot \mid
\cdot \right\rangle \right) \in Hilb\left( \mathbb{R}^{\left( a,b\right)
}\right) .$$
 Let $H$ be the reproducing kernel of \ $\mathcal{H}$ $\ .$We
have :
\begin{equation*}
\forall s,t\in \left( a,b\right) \ ,\ \text{\ }H\left( s,t\right)
=\sum_{j=0}^{n}\frac{\left( s-c\right) ^{j}}{j!}.\frac{\left( t-c\right) ^{j}%
}{j!}
\end{equation*}
\subsubsection{Example 2}
\qquad\\
 Let \ $a,b\in \mathbb{R}$ \ , \ $-\infty <a<b<+\infty $ \ and :
\\
$\mathcal{D}'(a,b)$ the space of distributions on $(a,b)$, 
\\
 $\mathcal{H}^{1}\left( a,b\right) =\left\{ x\in \mathcal{D}%
^{\prime }\left( a,b\right) \text{ ; \ }x\in L^{2}\left( a,b\right) \text{ \
and \ \ }x^{\prime }\in L^{2}\left( a,b\right) \right\} $ .
\\
 Let \ $\mathcal{H=}\left\{ x\in \mathcal{H}^{1}\left(
a,b\right) \text{ ; \ }x\left( a\right) =0\right\} $ \ and :
 $\forall x,y\in \mathcal{H}$ \ , \ $\left\langle x\mid
y\right\rangle =\int_{a}^{b}x^{\prime }\left( t\right) .y^{\prime }\left(
t\right) $ $dt$ . Then
$$\left( \text{\ }\mathcal{H}\text{ },\left\langle \cdot \mid
\cdot \right\rangle \right) \in Hilb\left( \mathbb{R}^{\left( a,b\right)
}\right).$$
 Let $H$ be the reproducing kernel of \ $\mathcal{H}.$
We have :
\begin{equation*}
\forall s,t\in \left( a,b\right) \ ,\ \ H\left( s,t\right) =\left(
t-s\right) _{+}+t-a=Min\left( t-a,s-a\right)
\end{equation*}
\subsubsection{Example 3 (Fourier kernel)}
\qquad\\
 Let \ $\mathcal{H=}\left\{ x\in \mathcal{H}^{1}\left(
0,2\pi \right) \text{ \ ; \ }x\left( 0\right) =x\left( 2\pi \right) \text{ \
and }\int_{0}^{2\pi }x\left( s\right) ds\text{\ }=0\right\} $ and

\qquad \qquad \qquad $\forall x,y\in \mathcal{H}$ \ , \ $\left\langle x\mid
y\right\rangle =\int_{0}^{2\pi }x^{\prime }\left( s\right) .y^{\prime
}\left( s\right) $ $ds$ . Then :
$$\left( \text{\ }\mathcal{H}\text{ },\left\langle \cdot \mid
\cdot \right\rangle \right) \in Hilb\left( \mathbb{R}^{\left( 0,2\pi \right)
}\right).$$
 Let $H$ be the reproducing kernel of \ $\mathcal{H}$.
 We
have :

\begin{equation*}
\forall s,t\in \left( 0,2\pi \right) \ ,\ \ H\left( s,t\right) =\left( \pi
\right) ^{-1}\left( \sum_{p=1}^{\infty }\frac{\cos p\left( t-s\right) }{\pi
p^{2}}\right)
\end{equation*}
\subsubsection{Example 4 ( the kernel is an operator)}
\qquad\\ Let \ $n\in \mathbb{N}^{\ast }$ , \ $\Omega $ \ an (non void)
open subset of \ $\mathbb{R}^{n}$ and \ $\alpha \in \mathbb{N}^{\ast }$ .
Let \ $\mathcal{H}^{m}\left( \Omega \right) $ \ be the space
of distributions which belong to \ $L^{2}\left( \Omega \right) $
 and which their derivatives of order less or equal to \ $%
\alpha $ \ belong to \ $L^{2}\left( \Omega \right) $.
 We set : 
\begin{equation*}
\forall x,y\in \mathcal{H}^{\alpha }\left( \Omega \right) \text{ , }\
\left\langle x\mid y\right\rangle _{\alpha }=\sum_{\left\vert p\right\vert
\leq \alpha }\left( \int_{\Omega }a_{p\text{ }}D^{p}x\left( \zeta \right)
.D^{p}y\left( \zeta \right) \text{ }d\zeta \right)
\end{equation*}
where \ $p=\left( p_{1},....,p_{n}\right) $ \ \ $\left\vert
p\right\vert =$\ $p_{1}+....+p_{n}$ . Then :
$$\left( \mathcal{H}^{\alpha }\left( \Omega
\right) ,\ \left\langle x\mid y\right\rangle _{\alpha }\right) \in \mathcal{D%
}^{\prime }\left( \Omega \right) .$$
 We denote by  $\mathcal{H}^{-\alpha }\left( \Omega \right) $%
 the topological dual \ of \ $\left( \mathcal{H}^{\alpha }\left(
\Omega \right) ,\ \left\langle x\mid y\right\rangle _{\alpha }\right) $
 which belongs to \ $\mathcal{D}^{\prime }\left( \Omega \right) 
$ .
 So, it is proved that the hilbertian kernel of \ $\mathcal{H}%
^{-\alpha }\left( \Omega \right) $\ (resp. $\mathcal{H}_{0}^{\alpha }\left(
\Omega \right) $ , $\mathcal{H}^{\alpha }\left( \Omega \right) $)
is the operator 
\begin{equation*}
D=\sum_{\left\vert p\right\vert \leq \alpha }\left( \left( -1\right) ^{\
\left\vert p\right\vert }a_{p\text{ }}D^{2p}\right)
\end{equation*}
(resp. its Green operator for the Dirichlet or Neumann
problem associated to \ $D$
 and the open set \ $\Omega $ .
\subsubsection{Example 5 ( the complex case)}
\qquad\\ Let \ $\Omega $ \ \ be a bounded (non void) open subset of \ $%
\mathbb{C}$ \ simply connected.
Given \ $m\in \mathbb{N}^{\ast }$, we denote by \ $\mathcal{A}%
_{m}\left( \Omega \right) $ \ the set of anlytical functions \ $f$ \ on \ $
\Omega$,
 such that :
\begin{equation*}
\int_{\Omega }\sum_{j=0}^{m}\left\vert f^{\left( j\right) }\left( z\right)
\right\vert ^{2}dz<+\infty
\end{equation*}
 where \ $z=x+iy$ \ \ and \ $%
\ $\ $\mathcal{A}_{0}\left( \Omega \right) =$\ \ $L^{2}\left( \Omega \right) 
$.
\begin{pp}
 We suppose that \ $\Omega $ \ is a disc centred at the origin
with radius $\rho $ .
\begin{itemize}
\item[(i)]
 We denote by \ $H_{0}$ \ the reproducing
kernel of \ \ $L^{2}\left( \Omega \right) $. 
Then :
\begin{equation*}
\forall t,z\in \Omega \text{ \ , \ }H_{0}\left( t,z\right) =\left( \pi \rho
^{2}\right) ^{-1}\left( 1-\frac{t\overline{z}}{\rho ^{2}}\right) ^{-1}
\end{equation*}
\item[(ii)]
 We denote by \ $H_{1}$ \ the reproducing
kernel of \ the set :
$\left\{ f\in \mathcal{A}_{1}\left( \Omega
\right) \ ;\ f\left( \zeta _{1}\right) =0\ \ \ \text{where }\ \zeta _{1}\in
\Omega \right\} $. Then :
\begin{equation*}
\forall t,z\in \Omega \text{ },\text{ }H_{1}\left( t,z\right) =-\left( \pi
\right) ^{-1}\left[ Log\left( 1-\frac{t\overline{z}}{\rho ^{2}}\right)
-Log\left( 1-\frac{t\overline{\zeta _{1}}}{\rho ^{2}}\right) -Log\left( 1-%
\frac{\zeta _{1}\overline{z}}{\rho ^{2}}\right) -Log\left( 1-\frac{\overline{%
\zeta _{1}\zeta _{1}}}{\rho ^{2}}\right) \right]
\end{equation*}
\end{itemize}
\end{pp}
\subsubsection{Tensor product of Hilbert spaces.}
\qquad\\
 Let \ $\left( E_{j},\left\langle \cdot \mid \cdot \right\rangle
_{j}\right) $, $j=1,.....,m$, $m\in \mathbb{N}^{\ast }$, be $m$ 
Hilbert spaces and
 $E~=~E_{1}\otimes~\ldots~\otimes~E_{m}$.
 We denote by \ $\left\langle \cdot \mid \cdot \right\rangle _{\otimes
}$ \ \ the scalar product on \ $E$ \ \ such that :

\qquad \qquad $\forall x_{j},y_{j}\in E_{j}$ \ , \ \ \ $j=1,.....,m$ , \ $%
\left\langle \otimes _{j=1}^{m}x_{j}\mid \otimes
_{j=1}^{m}y_{j}\right\rangle _{\otimes }$ $=\Pi _{j=1}^{m}\left\langle
x_{j}\mid y_{j}\right\rangle $.
\\
One can prove the
\begin{pp}
\quad\\
\begin{itemize}
\item[(i)] $\left( E,\left\langle \cdot \mid \cdot
\right\rangle _{\otimes }\right) $ \ is a (pre-)Hilbert space.
\item[(ii)] If \ $\left\{ e_{j}^{k}\text{ ; }k\in
J\left( j\right) \right\} $ is an (orthonormal) basis of $\left(
E_{j},\left\langle \cdot \mid \cdot \right\rangle _{j}\right) $,$j=1,œ\ldots, m$, then :
$$\qquad \qquad \left\{ e_{1}^{k_{1}}\otimes ....\otimes e_{m}^{k_{m}}%
\text{ ; }k_{1}\in J\left( 1\right) \text{ ,...., }k_{m}\in J\left( m\right)
\right\} $$ 
is an (orthonormal) basis of $\left( E,\left\langle \cdot
\mid \cdot \right\rangle _{\otimes }\right) $ .
\end{itemize}
\end{pp}
\subsection{(Hilbertian) Schwartz kernel of the tensor product of Hilbert
spaces.}
\subsubsection{The tensor norm \ $\sigma $}
\qquad\\ Let $\mathcal{A}_{j}$ \ be an e.l.c.s. which topology denoted \ $%
\mathcal{T}\left( \mathcal{A}_{j}\right) $ is defined by the filtering
family of (semi-)norms \ $\left\{ p_{j,k}\text{ ; }k\in J\left( j\right)
\right\} $ \ for \ $j=1,2$ .
\\
 We call projective (resp. inductive ) tensor product of \ $\left( 
\mathcal{A}_{1}\text{ ,}\mathcal{T}\left( \mathcal{A}_{1}\right) \right) $ \
by \ $\left( \mathcal{A}_{2}\text{ ,}\mathcal{T}\left( \mathcal{A}%
_{2}\right) \right) $
the e.l.c.s.
$$\mathcal{A}_{1}\pi \mathcal{A}_{2}\underset{def}{=}\left( 
\mathcal{A}_{1}\otimes \mathcal{A}_{2}\text{ , \ }\left\{ \left( p_{k,1}
\text{ }\pi \text{ }p_{l,2}\right) \text{ , }\left( k,l\right) \in J\left(
1\right) \times J\left( 2\right) \right\} \right) $$
$$\left( resp.\quad\mathcal{A}_{1}\varepsilon \mathcal{A}_{2}\underset{def}{=}%
\left( \mathcal{A}_{1}\otimes \mathcal{A}_{2}\text{ , \ }\left\{ \left(
p_{k,1}\text{ }\varepsilon \text{ }p_{l,2}\right) \text{ , }\left(
k,l\right) \in J\left( 1\right) \times J\left( 2\right) \right\} \right)
\right) $$
We shall denote by \ $\mathcal{A}_{1}\widehat{\pi }\mathcal{A}_{2}$ \ 
$\left( resp.\mathcal{A}_{1}\widehat{\varepsilon }\mathcal{A}_{2}\right) $ \
a completion of \ \ $\mathcal{A}_{1}\pi \mathcal{A}_{2}$ \ $\left( resp.%
\mathcal{A}_{1}\varepsilon \mathcal{A}_{2}\right) $.
It can be proved ( cf \ M.A....) that there exists on \ $\mathcal{A}%
_{1}\otimes \mathcal{A}_{2}$ \ a reasonable tensor norm $\sigma $
which is the arithmetico-geometric of \ $\varepsilon $ \ and \ $\pi $ .
\\
 We shall denote by \ $\mathcal{A}_{1}\sigma \mathcal{A}_{2}$ \ the
linear space \ $\mathcal{A}_{1}\otimes \mathcal{A}_{2}$ \ with the tensor
norm $\sigma $ \ and
by \ $\mathcal{A}_{1}\widehat{\sigma }\mathcal{A}_{2}$ \ a completion of \ $%
\mathcal{A}_{1}\sigma \mathcal{A}_{2}$ .
\\
So, we can easily prove the following injections :
$$\mathcal{A}_{1}\widehat{\pi }\mathcal{A}%
_{2}\rightarrow \mathcal{A}_{1}\widehat{\sigma }\mathcal{A}_{2}\rightarrow 
\mathcal{A}_{1}\widehat{\varepsilon }\mathcal{A}_{2}.$$
That property can be extended to the tensor product \ $\mathcal{A}%
_{1}\otimes \mathcal{A}_{2}\otimes ....\otimes \mathcal{A}_{m}$
 $m$ \ e.l.c.s., \ $m\in \mathbb{N}^{\ast }$ , \ $m\geq 3$ .
\subsection{(Hilbertian) Schwartz kernel of a tensor product}
\qquad \\ Let \ $\left( E_{j},\left\langle \cdot \mid \cdot
\right\rangle _{j}\right) $ \ \ $j=1,2$ \ two Hilbert spaces (on \ $\mathbb{K%
}$ ) and \ $E=E_{1}\otimes E_{2}$.
 From the above paragraph we can deduce that : \ $\left(
E,\left\langle \cdot \mid \cdot \right\rangle _{\otimes }\right)
=E_{1}\sigma E_{2}$.
Moreover, if \ \ $\left( E_{j},\left\langle \cdot \mid \cdot
\right\rangle _{j}\right) \in Hilb\left( \mathcal{A}_{j}\right) $\ $\ $, \ $%
j=1,2$ \ \ then: \ \ $\left( E,\left\langle \cdot \mid \cdot \right\rangle
_{\otimes }\right) \in Hilb\left( \mathcal{A}_{1}\widehat{\sigma }\mathcal{A}%
_{2}\right) $.
 Thus, if \ \ $\mathcal{E}_{j}$ \ is the Schwartz kernel of \ \ 
$\left( E_{j},\left\langle \cdot \mid \cdot \right\rangle _{j}\right) \ \ ,\
j=1,2$ , then :
\\
 the kernel of $\ \ \left( E,\left\langle \cdot \mid \cdot
\right\rangle _{\otimes }\right) $ \ \ is equal \ to \ 
$\mathcal{E}_{1}\otimes \mathcal{E}_{2}$ 
 and
$$\forall a_{k}^{\ast },b_{k}^{\ast }\in \mathcal{A}_{k}^{\ast },
\quad k=1,2,\quad
 \left\langle \mathit{\ }\left( \mathcal{E}_{1}\otimes 
\mathcal{E}_{2}\right) \left( a_{1}^{\ast }\otimes a_{2}^{\ast }\right)
,\left( b_{1}^{\ast }\otimes b_{2}^{\ast }\right) \right\rangle
_{\circledcirc }=\left\langle \mathit{\ }\mathcal{E}_{1}a_{1}^{\ast }\text{ }%
,b_{1}^{\ast }\right\rangle _{1}.\left\langle \mathit{\ }\mathcal{E}%
_{2}a_{2}^{\ast }\text{ },b_{2}^{\ast }\right\rangle_{2} .
$$
 That last property can be extended to the tensor product of \ the $m$
\ hilbertian subspaces
$\left( E_{j},\left\langle \cdot \mid \cdot \right\rangle _{j}\right) 
$ \ \ of \ the e.l.c.s. $\mathcal{A}_{j}$ \ for \ \ $j=1,\ldots ,m$ .
\subsection{Hilbertian Fock spaces and their Schwartz kernels.}
\subsubsection{Hilbertian Fock spaces}
\qquad \\
Using the same notations as those in the paragraph B4, we set :
 Let \ $\left( E,\left\langle \cdot \mid \cdot \right\rangle \right) $
\ be a Hilbert space \ (on \ $\mathbb{K}$ ) \ and \ $\left\Vert \cdot
\right\Vert =\sqrt{\left\langle \cdot \mid \cdot \right\rangle }$.
\\
For $\ \ m\in \mathbb{N}^{\ast }$ , \ $m\geq 2$ .We denote by\ \ \ \ $\left(
E^{\vee m},\left\langle \cdot \mid \cdot \right\rangle _{\vee }\right) $ \ \ 
$(resp.$\ \ $\left( E^{\wedge m},\left\langle \cdot \mid \cdot \right\rangle_{\wedge }\right) $\ , the (pre-) Hilbert space
such that : \ \ $\forall x_{j},y_{j}\in E$,\ $j=1,\ldots,m$,
$$
\left\langle \vee _{j=1}^{m}x_{j}\mid \vee
_{j=1}^{m}y_{j}\right\rangle _{\vee }=
\frac{1}{m!}
\left(
\sum_{\rho \in G_{m}}\text{ }\left\langle x_{1}\mid y_{\rho \left( 1\right)
}\right\rangle \ldots\left\langle x_{m}\mid y_{\rho \left( m\right)
}\right\rangle \right) $$
$$\left (resp.\quad
\left\langle \wedge _{j=1}^{m}x_{j}\mid \wedge
_{j=1}^{m}y_{j}\right\rangle _{\vee }=
 \frac{1}{m!}
\left(
\sum_{\rho \in G_{m}}\text{ }\varepsilon \left( \rho \right) \text{ }%
\left\langle x_{1}\mid y_{\rho \left( 1\right) }\right\rangle
......\left\langle x_{m}\mid y_{\rho \left( m\right) }\right\rangle \right)
\right).$$
\subsubsection{Schwartz kernels of hilbertian Fock spaces}
\qquad\\
 Let \ $\mathcal{A}$ \ be an e.l.c.s. \ . We suppose that \ $\left(
E,\left\langle \cdot \mid \cdot \right\rangle \right) \in Hilb\left( 
\mathcal{A}\right) $ \ and that \ $\mathcal{E}$ \ is
its Schwartz kernel. So, we denote by \ $\mathcal{E}^{\vee m}$ \ and \ \ $\mathcal{E%
}^{\wedge m}$\ \ the Schwartz kernels \ of \ 
$\left( E^{\vee m},\left\langle \cdot \mid \cdot \right\rangle _{\vee
}\right) $ \ and $E^{\wedge m},\left\langle \cdot \mid \cdot \right\rangle_{\wedge }$ \ respectively.
\\
We know that : 
if \ \ $\mathcal{A}_{\sigma }^{\otimes m}=\mathcal{A\sigma A\sigma }%
......\sigma \mathcal{A}$ \ where \ $\sigma $ \ appears \ $\left( m-1\right) 
$ times, then :
\\
$\left( E^{\vee m},\left\langle \cdot \mid \cdot \right\rangle _{\vee
}\right) $ \ and $E^{\wedge m},\left\langle \cdot \mid \cdot \right\rangle
_{\wedge }$ \ belong to \ $Hilb\left( \mathcal{A}_{\sigma }^{\otimes
m}\right) $ \ and :
\\
 $\forall a_{k}^{\ast },b_{k}^{\ast }\in \mathcal{A}^{\ast }$ \ , $\
k=1,\ldots,m$,
$$\left\langle \ \mathcal{E}^{\vee m}\left( a_{1}^{\ast }\vee
\ldots\vee a_{m}^{\ast }\right) ,\left( b_{1}^{\ast }\vee \ldots
\vee
b_{m}^{\ast }\right) \right\rangle _{\vee }=\frac{1}{m!}\left(
\sum_{\rho \in G_{m}}\text{ }\left\langle \mathcal{E}a_{1}^{\ast }\mid
b_{\rho \left( 1\right) }^{\ast }\right\rangle 
\ldots\left\langle \mathcal{E}
a_{m}^{\ast }\mid b_{\rho \left( m\right) }^{\ast }\right\rangle \right) $$
and
$$\left\langle \ \mathcal{E}^{\wedge m}\left( a_{1}^{\ast }\vee
\ldots
\vee a_{m}^{\ast }\right) ,\left( b_{1}^{\ast }\vee 
\ldots
\vee
b_{m}^{\ast }\right) \right\rangle_{\vee }
=
\frac{1}{m!}
\left(
\det \left[ a_{j}^{\ast },b_{k}^{\ast }\right] _{1\leq j,k\leq m}\right) .$$
The following lemma can be proved easily :
\begin{lm}
 Let \ \ $\mathcal{W}_{j}$ \ be an e.l.c.s, $\ \ \ \left(
F_{j},\left\langle \cdot \mid \cdot \right\rangle _{j}\right) $ \ an
hilbertian subspace of $\left( \mathcal{W}_{j}\right) $ \ \ \ \ \ \ 
 and \ $\Phi _{j}$ \ its Schwartz kernel for \ \ $j=1,\ldots,m$.
Let \ $\left( F,\left( \mid \right) \right) $ \ be the cartesian
product~:  $\left( F_{1},\left\langle \cdot \mid \cdot \right\rangle
_{1}\right) \times ......\times \left( F_{m},\left\langle \cdot \mid \cdot
\right\rangle _{m}\right) $ .
Then $\ $: $\ \ \left( F,\left( \mid \right) \right) \in Hilb\left( 
\mathcal{W}_{1}\times .......\times \mathcal{W}_{m}\right) $\ \ and if \ $%
\Phi $ \ is its Schwartz kernel
we have :
$$
\Phi =
\left[ 
\begin{array}{ccccc}
\Phi _{1} &  & 0 & 0 & 0 \\ 
0 & \Phi _{2} & 0 & 0 & 0 \\ 
0 & 0 & ... & 0 & 0 \\ 
... & ... & ... & ... & ... \\ 
0 & 0 & 0 & 0 & \Phi _{m}%
\end{array}
\right]. $$
\end{lm}
 So, with the same notations as at the beginning of this paragrah we
deduce that :
\\
if \ \ $\mathcal{F}_{\gamma }^{\bot }\left( E,\left\Vert \cdot
\right\Vert \right) $ \ \ is the (tensorial) Fock space associated to \ $%
\left( E,\left\Vert \cdot \right\Vert \right) $ \ with the r.t.n. $\gamma $
 and \ \ $\Gamma_{\gamma }^{\bot }\left( 
\mathcal{E}\ \right) $
its Schwartz kernel then :

 $\Gamma_{\gamma }^{\bot }\left( \mathcal{E}\
\right) $  is matrix diagonal block-matrix and its diagonal is equal to 
$1_{\mathbb{R}}$ $\mathcal{E}$ $\mathcal{E}^{\bot 2}\ldots
\mathcal{E}^{\bot m}.$
\\
 \textbf{Exercise :}
\\
\qquad \qquad Calculate \ $\Gamma_{\gamma }^{\bot }\left( 
\mathcal{E}\ \right) $ \ when \ $\mathcal{E}$ \ is one of the five examples
of Schwartz kernels in the
paragraph 4 above.
\section*{Bibliography}

$[1] $ \textsc{P. Angles, }The structure of the Clifford algebra.

\qquad\ Advanced applied Clifford algebra, vol.19,
n$^\circ$3-4, p.585-610, 2009.

$[2] $ \textsc{M.\ Atteia,} Hilbertian kernels and spline
functions.

\qquad\ Studies of computational Mathematics 4, North Holland,

$[3] $ \textsc{R.\ Deheuvels,} Formes quadratiques et groupes
classiques.

\qquad\ Presses universitaires de France, 1981.

$[4] $ \textsc{A.\ Grothendieck,} Produits tensoriels
topologiques et espaces nucl\'{e}aires.

\qquad AMS 16, 1955.

$[5] $ \textsc{P. Lounesto,} Clifford algebra and spinors.

\qquad Cambridge University Press, 1997.

$[6] $ J Neveu, Processus al\'{e}atoires gaussiens.

\qquad Presses de l'Universit\'{e} de Montr\'{e}al, 
n$^\circ$34, 1968.

$[7] $ \textsc{R.\ Ryan,} Clifford algebra in Analysis and
related topics.

\qquad Studies in Advanced Mathematics, CRC press, 1996.

\pagebreak

$[8] $ \textsc{R.\ Ryan,} Introduction to tensor products of
Banach space.
 
\hspace{1cm} Springer, 2002.

$[9]$ \textsc{L.\ Schwartz,} Sous-espaces hilbertiens d'espaces
vectoriels topologiques

\qquad et noyaux associ\'{e}s.

\qquad Journal d'Analyse Math\'{e}matique, J\'{e}rusalem, vol.13, 1964.\ 
\end{document}